\documentclass[11pt]{article}
\usepackage{color}
\usepackage{tikz}
\usepackage{amssymb, amsmath}
\usepackage{epsfig}
\usepackage{graphicx}

\newtheorem{example}{Example}

\begin{document}

\title{On the computation of the nth power of a matrix}

\author{Nikolaos Halidias \\
{\small\textsl{Department of Mathematics }}\\
{\small\textsl{University of the Aegean }}\\
{\small\textsl{Karlovassi  83200  Samos, Greece} }\\
{\small\textsl{email: nikoshalidias@hotmail.com}}}

\maketitle

\begin{abstract} In this note we discuss the problem of finding
the nth power of a matrix which is strongly connected to the study
of Markov chains and others mathematical topics. We observe the
known fact (but maybe not well known) that the Cayley-Hamilton
theorem is of key importance to this goal. We also demonstrate the
classical Gauss elimination technique as a tool to compute the
minimum polynomial of a matrix without necessarily know  the
characteristic polynomial.
\end{abstract}


{\bf 2010 Mathematics Subject Classification}

\section{Using Cayley-Hamilton theorem to find the nth power of a matrix}
Let  $A_{m \times m}$ be a  real matrix and suppose that we want
to compute the nth power. Let us denote by $\delta(k)$ the
characteristic polynomial of $A$. Then it holds that
\begin{eqnarray}\label{longdivision}
k^n = \delta_{A}(k) \pi (k) + v(k)
\end{eqnarray}
 where  $v(k)$ is a polynomial of degree less or equal to
$m-1$. Using the Cayley-Hamilton theorem we see that $A^n = v(A)$
since $\delta(A) = 0$. Therefore, in order to find the form of
$A^n$ we have to find the coefficients of $v(k)$ which are at most
$m$. In this direction we will use the eigenvalues
$k_1,\cdots,k_n$ (not necessarily distinct) of the matrix $A$
since $\delta(k_i) = 0$. Setting $k = k_i$ we produce some
equations involving the unknown coefficients of $v$. Of course if
some of the eigenvalues are of multiplicity two or more then we
have to produce some more equations. This can be done by
differentiating equation \ref{longdivision} and setting $k = k_i$.
We differentiate this equation $m_i-1$ times (where $m_i$ is the
multiplicity of $k_i$ eigenvalue) and each time we set $k = k_i$
in order to produce one more equation. We do this for any
eigenvalue $k_j$ with multiplicity $m_j$. In this fashion we
produce $m$ different equations in order to determine the $m$
unknown coefficients of $v$. One can easily observe that the
computation of the coefficients of $v(k)$ is in fact an
interpolation problem which has a unique solution (see \cite{Dav},
Thm. 2.2.2). Therefore, we can always find the nth power of matrix
by using the Cayley-Hamilton theorem even if the matrix is not
diagonizable. Below, we give some examples.

\begin{example} We will compute the nth power of the matrix
\begin{eqnarray*}
P = \left(%
\begin{array}{rrr}
  -3 & 6 & 0 \\
  2 & 1 & 0 \\
  0 & 0 & 3 \\
\end{array}%
\right)
\end{eqnarray*}
The eigenvalues are
 $k_1 = 3$ with multiplicity two and  $k_2 = -5$ that is the characteristic polynomial is  $\delta_P(k) = (k-3)^2(k+5)$.

Therefore, it holds that  $k^n = t(k) \delta_P(k) + v(k)$ where
$v(k) = a k^2 + bk +c$. We have to produce three equations in
order to determine the unknown coefficients. Setting $k_1 = 3$ we
obtain the first equation $3^n = 9 a + 3 b +c$ and setting  $k=-5$
we produce another one which is   $(-1)^n 5 = 25 a - 5 b + c$. In
order to have one more equation we differentiate the equality $k^n
= t(k) \delta_P(k) + v(k)$  and then we set again $k=3$ because
this eigenvalue has multiplicity two. Therefore, we obtain the
third equation which is
 $3^{n-1} n  = 6a + b$. Next we solve for $a,b,c$ and thus
\begin{eqnarray*}
P^n = a P^2 + b P + c = \left(%
\begin{array}{ccc}
  21a-3a+c & -12a+6b & 0 \\
  -4a+2b & 13a+b+c & 0 \\
  0 & 0 & 9a + 3b + c \\
\end{array}%
\right)
\end{eqnarray*}
Finally, one can verify the result by induction.
\end{example}

\begin{example}\label{Cayleymigadikoi} We will compute the nth
power of the following matrix
\begin{eqnarray*}
P = \left(%
\begin{array}{ccc}
  0 & 1 & 0 \\[0.2cm]
  0 & \frac{2}{3} & \frac{1}{3} \\[0.2cm]
  \frac{1}{3} & 0 & \frac{2}{3} \\
\end{array}%
\right)
\end{eqnarray*}
This matrix has the following eigenvalues,
 $k_1 = 1$, $k_2 = \frac{1 + i
\sqrt{3}}{6}$ and $k_3 = \frac{1 - i \sqrt{3}}{6}$. We have the
following equation  $k^n = \delta_P(k) t(k) + v(k)$ where  $v(k) =
a k^2 + bk + c$. Next, we will produce three equations in order to
evaluate the unknown coefficients. We transform every complex
number in the polar form, that is
 $x+iy = r(\cos \phi + i \sin \phi)$ where  $r =
\sqrt{x^2+y^2}$ and $\phi = \arctan \frac{y}{x}$. Therefore we
have that
\begin{eqnarray*}
k_2 & = & \frac{1}{3} \left( \cos \frac{\pi}{3} + i \sin
\frac{\pi}{3}\right)
\\
k_3 & =& \frac{1}{3} \left( \cos \frac{\pi}{3} - i \sin
\frac{\pi}{3}\right)
\end{eqnarray*}

We obtain the following equations by setting $k = k_i$ in the
equation $k^n = \delta_P(k) t(k) + v(k)$
\begin{eqnarray*}
1^n & = & a + b+ c \\[0.2cm]
 \frac{1}{3^n} \left( \cos \frac{n\pi}{3} + i \sin
\frac{n\pi}{3}\right) & = & a \frac{1}{3^2} \left( \cos
\frac{2\pi}{3} + i \sin \frac{2\pi}{3}\right) + b \frac{1}{3}
\left( \cos \frac{\pi}{3} + i \sin \frac{\pi}{3}\right) + c
\\[0.2cm]
\frac{1}{3^n} \left( \cos \frac{n\pi}{3} - i \sin
\frac{n\pi}{3}\right) & = & a \frac{1}{3^2} \left( \cos
\frac{2\pi}{3} - i \sin \frac{2\pi}{3}\right) + b \frac{1}{3}
\left( \cos \frac{\pi}{3} - i \sin \frac{\pi}{3}\right) + c
\end{eqnarray*}
Adding the second equation to the third we produce the equality
\begin{eqnarray*}
\frac{1}{3^n} \cos \frac{n \pi}{3} = a \frac{1}{3^2}\cos
\frac{2\pi}{3} + b \frac{1}{3}\cos \frac{\pi}{3} + c
\end{eqnarray*}
while subtracting the third equation from the second we produce
the following
\begin{eqnarray*}
\frac{1}{3^n} \sin \frac{n \pi}{3} = a \frac{1}{3^2}\sin
\frac{2\pi}{3} + b \frac{1}{3}\sin \frac{\pi}{3}
\end{eqnarray*}
Finally, we have the following system of equations
\begin{eqnarray*}
1^n & = & a + b+ c \\[0.2cm]
\frac{1}{3^n} \cos \frac{n \pi}{3} &  = & a \frac{1}{3^2}\cos
\frac{2\pi}{3} + b \frac{1}{3}\cos \frac{\pi}{3} + c\\[0.2cm]
\frac{1}{3^n} \sin \frac{n \pi}{3} & = & a \frac{1}{3^2}\sin
\frac{2\pi}{3} + b \frac{1}{3}\sin \frac{\pi}{3}
\end{eqnarray*}
Solving for $a,b,c$  we obatin
\begin{eqnarray*}
a & = & - \frac{3 \sqrt{3}}{7} \left( \sqrt{3} \frac{1}{3^n} \cos
\frac{n \pi}{3} - \sqrt{3} + 5 \frac{1}{3^n} \sin \frac{n \pi}{3}
\right) \\[0.2cm]
b & = &  \frac{ \sqrt{3}}{7} \left( \sqrt{3} \frac{1}{3^n} \cos
\frac{n \pi}{3} - \sqrt{3} + 19 \frac{1}{3^n} \sin \frac{n \pi}{3}
\right) \\[0.2cm]
c & = &  \frac{\sqrt{3}}{21} \left( 6 \sqrt{3} \frac{1}{3^n} \cos
\frac{n \pi}{3} + \sqrt{3} -12  \frac{1}{3^n} \sin \frac{n \pi}{3}
\right)
\end{eqnarray*}
 and therefore $P^n = a P^2 + b P +c I_{2 \times 2}$. One can
 verify the result by induction.
\end{example}

\section{Minimum polynomial} Recall that the minimum polynomial of
$A$ is the polynomial $q(k)$ of less degree such that $q(A) = 0$.
Therefore one can use the minimum polynomial rather than the
characteristic one in order to compute the nth power of the matrix
$A$.

We will demonstrate the classical Gauss elimination procedure in
order to find  the minimum polynomial.

Let the minimum polynomial has the form
\begin{eqnarray*}
q(k) = k^{r} + a_{r-1} k^{r-1} + \cdots + a_0
\end{eqnarray*}
ме $r \leq m$.

Obviously the relation
\begin{eqnarray}\label{sistimaeksisoseongrammikosane}
b_{r-1}A^{r-1} + \cdots + b_0 I_{m \times m} = 0_{m \times m}
\end{eqnarray}
drive us to the conclusion that  $b_{r-1} = b_{r-2} = \cdots = b_0
=0$ otherwise there will be a polynomial $\hat{q}$ such that
$\hat{q}(A) = 0$ with less degree than the minimum polynomial
which is a contradiction. That means that the matrices
 $I_{m \times m}, A,
\cdots, A^{r-1}$ are linearly independent while the matrices $I_{n
\times n}, A, \cdots, A^{r-1}, A^r$ are linearly dependent.
Obviously the matrices
$$I_{m \times m}, A, \cdots, A^{r-1}, A^r, A^{r+1}, \cdots,A^m$$
are also linearly dependent.

 Let  $C_{m^2 \times r}$ the matrix that has in its first column the identity  matrix $I$, that is
 at the first
$m$ places of the first column of $C$ we put the first column of
$I$, next at the next $m$ places of the first column of $C$ we put
the second column of $I$ and so on. We do the same for the
matrices $A,A^2, \cdots, A^{r-1}$. This matrix is the matrix of
the system \ref{sistimaeksisoseongrammikosane}. Since the matrices
 $I,A,A^2,\cdots,A^{r-1}$ are linearly dependent then the reduced row echelon form of
 $C$ will have the following form
\begin{eqnarray*}
\hat{C} = \left(%
\begin{array}{cccc}
  1 & 0 & 0 & \cdots \\
  0 & 1 & 0 & \cdots \\
  0 & 0 & 1 & \cdots \\
  \vdots & \vdots & \vdots & \vdots \\
 0 & 0 & 0 & 0 \\
  0 & 0 & 0 & 0 \\
\end{array}%
\right)
\end{eqnarray*}
where the number of leading 1 equals to
 $r$.

Let now   $D_{m^2 \times (m+1-r)}$ the matrix that has in its
columns the matrices  п $A^r,\cdots,A^m$ with the same fashion as
with $C$. We construct next the matrix $B = (C|D)$ (which is in
fact the matrix of system \ref{ejisoseiselaxistoupolionimou}
below) and we compute the reduced row echelon form. Obviously, the
number of the leading 1 equals $r$ again and these are located at
the $r$ first columns of $B$.

That means that the system
\begin{eqnarray}\label{ejisoseiselaxistoupolionimou}
a_n A^n + a_{n-1} A^{n-1} + \cdots + a_0I_{n \times n} = 0_{n
\times n}
\end{eqnarray}
has infinite many solutions with  $n+1-r$ free parameters.
 Setting $a_r = 1$,
$a_{r+1} = a_{r+2} = \cdots = a_{n} = 0$ and solving for the
others we obtain the minimum polynomial
\begin{eqnarray*}
q(k) = a_0 + a_1 k + \cdots + a_{r-1} k^{r-1} + k^r
\end{eqnarray*}
One can easily verify that the minimum polynomial is as follows
\begin{eqnarray*}
q(k) = k^{r} - \hat{B}_{r,r+1} k^{r-1} - \hat{B}_{r-1,r+1} k^{r-2}
- \cdots - \hat{B}_{1,r+1}
\end{eqnarray*}

\begin{example} Let the matrix
\begin{eqnarray*}
A = \left(%
\begin{array}{rrr}
  -4 & 2 & 0 \\
  -2 & -1 & 0 \\
  0 & 0 & 1 \\
\end{array}%
\right)
\end{eqnarray*}
We will compute the minimum polynomial.

We construct the matrix $B$ by using the matrices $I,A,A^2,A^3$.
Then
\begin{eqnarray*}
B= \left(%
\begin{array}{rrrrr}
  1 & -4 & 12 & -28 \\
  0 & -2 & 10 & -34 \\
  0 & 0 & 0 & 0 \\
  0 & 2 & -10 & 34 \\
  1 & -1 & -3 & 23 \\
  0 & 0 & 0 & 0 \\
  0 & 0 & 0 & 0 \\
  0 & 0 & 0 & 0 \\
  1 & 1 & 1 & 1 \\
\end{array}%
\right)
\end{eqnarray*}
end the reduced row echelon matrix is the following
\begin{eqnarray*}
\hat{B} = \left(%
\begin{array}{rrrr}
  1 & 0 & 0 & 8 \\
  0 & 1 & 0 & -3 \\
  0 & 0 & 1& -4 \\
  0 & 0 & 0 & 0 \\
  \vdots & \vdots & \vdots & \vdots \\
\end{array}%
\right)
\end{eqnarray*}
The number of the leading 1 is 3 and that means that the minimum
polynomial if of third degree. Thus,
\begin{eqnarray*}
q(k) = \delta(k) = k^3 + 4k^2 +3k - 8
\end{eqnarray*}
and in this case coincides with the characteristic polynomial.
\end{example}

\begin{example} Let the matrix
\begin{eqnarray*}
A = \left(%
\begin{array}{rrr}
  -3 & 6 & 0 \\
  2 & 1 & 0 \\
  0 & 0 & 3 \\
\end{array}%
\right)
\end{eqnarray*}
We construct the matrix $B$ as before and therefore we have
\begin{eqnarray*}
B= \left(%
\begin{array}{rrrrr}
  1 & -3 & 21 & -87 \\
  0 & 2 & -4 & 38 \\
  0 & 0 & 0 & 0 \\
  0 & 6 & -12 & 114 \\
  1 & 1 & 13 & -11 \\
  0 & 0 & 0 & 0 \\
  0 & 0 & 0 & 0 \\
  0 & 0 & 0 & 0 \\
  1 & 3 & 9 & 27 \\
\end{array}%
\right)
\end{eqnarray*}
The reduced row echelon matrix of $B$ is
\begin{eqnarray*}
\hat{B} = \left(%
\begin{array}{rrrr}
  1 & 0 & 15 & -30 \\
  0 & 1 & -2 & 19 \\
  0 & 0 & 0& 0 \\
  \vdots & \vdots & \vdots & \vdots \\
\end{array}%
\right)
\end{eqnarray*}
The number of the leading 1 is 2 therefore the degree of the
minimum polynomial is 2. Setting
  $a_2 = 1$, $a_3 = 0$ and solving for the other we get
\begin{eqnarray*}
q(k) = k^2 +2k -15 = (k+5)(k-3)
\end{eqnarray*}
That means that the matrix  $A$ has the  $k_1 = -5$ and  $k_2 = 3$
as eigenvalues.
\end{example}

\end{document}